\documentclass[conference]{IEEEtran}

\usepackage{amsmath,amssymb,amsthm}
\usepackage{graphicx}
\usepackage{bm,color}
\usepackage{url}
\usepackage{cite}
\usepackage{textcomp}
\usepackage{slashbox}

\theoremstyle{definition}

\IEEEoverridecommandlockouts

\begin{document}
\title{Risk-Aware Management\\ of Distributed Energy Resources}

\author{\IEEEauthorblockN{Yu Zhang, Nikolaos Gatsis, Vassilis Kekatos, and Georgios B. Giannakis}
\authorblockA{Dept. of ECE and DTC, University of Minnesota, Minneapolis, USA \\
Emails: \{yuzhang, gatsisn, kekatos, georgios\}@umn.edu}
\thanks{This work was supported by the Inst. of Renewable Energy and the Environment
(IREE) grant no. RL-0010-13, Univ. of Minnesota.}}

\maketitle

\begin{abstract}
High wind energy penetration critically challenges the economic dispatch of current and future power systems. Supply and demand must be balanced at every bus of the grid, while respecting transmission line ratings and accounting for the stochastic nature of renewable energy sources. Aligned to that goal, a network-constrained economic dispatch is developed in this paper. To account for the uncertainty of renewable energy forecasts, wind farm schedules are determined so that they can be delivered over the transmission network with a prescribed probability. Given that the distribution of wind power forecasts is rarely known, and/or uncertainties may yield non-convex feasible sets for the power schedules, a scenario approximation technique using Monte Carlo sampling is pursued. Upon utilizing the structure of the DC optimal power flow (OPF), a distribution-free convex problem formulation is derived whose complexity scales well with the wind forecast sample size. The efficacy of this novel approach is evaluated over the IEEE 30-bus power grid benchmark after including real operation data from seven wind farms.
\end{abstract}

\section{Introduction}
The scarcity and the environmental impact of conventional energy resources raise major concerns worldwide, and drive industry to aggressively incorporate renewable energy, which is sustainable and clean. Coming from natural resources such as wind, sunlight, biomass, and geothermal heat, renewable energy-based electricity production has been developing rapidly in the past decade. Wind power generation for instance is growing at an annual rate of $20$\%, and has already met a worldwide installed capacity of $282.5$ GW by the end of 2012~\cite{GWEC}. The U.S. Department of Energy proposed and examined a goal of using wind energy to generate $20$\% of the nation's electricity demand by 2030~\cite{DOE08}.

The goal of high renewable energy penetration is challenged by the stochastic availability and intermittency of renewable energy, which must be accounted for by system operators during scheduling of generation, reserves, and dispatchable loads. This paper develops a chance-constrained optimization approach for economic power scheduling with renewables, with focus on controlling the risk stemming from potentially inadequate supply of renewable energy.

Prior works have dealt with the supply-demand imbalance issue under the uncertain supply of renewables. Single-period, chance-constrained economic dispatch is studied for a power system with both thermal generators and wind turbines in~\cite{LiuX10}. By using a here-and-now approach, a loss-of-load probability (LOLP)-guaranteed dispatch is obtained. A stochastic programming approach for economic dispatch simultaneously penalizing overestimation and underestimation of wind power is investigated in~\cite{HetzerYB08}. Multi-period economic dispatch with spatiotemporal wind forecasts is pursued in~\cite{XieGZG11}. By upper bounding wind power schedules by their forecasts, a deterministic optimization formulation is derived. Its solution though can be very sensitive to the accuracy of the wind power forecast.
Chance-constrained multi-period economic dispatch with multiple correlated wind farms has been explored recently in~\cite{YuNGGG-ISGT13}.

All aforementioned works limit their focus on the economic dispatch problem, which ignores the transmission network. Accounting for the transmission network leads to the optimal power flow (OPF) problem, which includes balance constraints for every network node, and flow limit constraints for every line; see e.g.,~\cite{ChWoWa00}. If not properly considered during system scheduling, the aforementioned constraints are more likely to be violated due to the stochastic and intermittent nature of renewable energy injections. To this end, relying upon Gaussianity assumptions for the wind power output and conic programming techniques, chance-constrained optimal power flow has been recently pursued in~\cite{BiChHa12} and~\cite{SjGaTo12}.

This paper deals with chance-constrained DC OPF for power systems with renewables; but different from~\cite{BiChHa12} and~\cite{SjGaTo12}, it develops a scheduling methodology that does not rely on Gaussianity. To address the stochastic nature of renewable energy, the proposed formulation introduces scheduled renewable energy injections as design variables, and allows the actual harvested energy to be inadequate with low risk. To effectively cope with the intractability of risk constraints, the proposed algorithm builds upon the scenario approximation of~\cite{YuNGGG-ISGT13}. Numerical tests are performed to corroborate the effectiveness of the novel approach using real wind farm operation data, and the IEEE 30-bus power grid benchmark \cite{kaggle}, \cite{PSTCA}, \cite{MATPOWER}.


The remainder of the paper is organized as follows. Section~\ref{sec:Probformulation} formulates the risk-constrained energy management problem, followed by the development of the scenario approximation approach in Section~\ref{Section:SAA}. Numerical results are reported in Section~\ref{sec:Numericalresults}, and conclusions are drawn in Section~\ref{sec:Conclusion}.

\section{Risk-Aware Energy Management}\label{sec:Probformulation}
Consider a power system with $M$ buses. Let $p_{G_m}$ denote the power output of a thermal generator  and $p_{D_m}$ the power dissipation of a load, both residing at bus $m$. While $p_{G_m}$ is a decision variable, load $p_{D_m}$ is considered fixed here for simplicity. Due to plant limitations, the generator power output is constrained to lie between lower and upper bounds $p_{G_m}^{\min}$ and $p_{G_m}^{\max}$, respectively. Furthermore, if a renewable energy producer is located at bus $m$, two quantities will be associated with it: the predicted wind power generation $z_m$, and the power $w_m$ scheduled to be injected to bus $m$. Note that the former is a random variable, whereas the latter is a decision variable. Define further the $M$-dimensional vectors
\begin{align*}
\mathbf{p}_G & :=[p_{G_1}~\ldots~p_{G_M}]^T \\
\mathbf{p}_D & :=[p_{D_1}~\ldots~ p_{D_M}]^T \\
\mathbf{z} & := [z_1~\ldots~z_M]^T \\
\mathbf{w} & := [w_1~\ldots~ w_M]^T \\
\mathbf{p}_{G}^{\min} & := [p_{G_1}^{\min}~\ldots~p_{G_M}^{\min}]^T \\
\mathbf{p}_{G}^{\max} & := [p_{G_1}^{\max}~\ldots~p_{G_M}^{\max}]^T
\end{align*}
where $(\cdot)^T$ denotes transposition. With these definitions, the nodal injections into the transmission grid can be expressed in vector form as $\mathbf{p}_G + \mathbf{w} - \mathbf{p}_D$.

Focusing next on the transmission network, let $L$ denote the number of lines in the grid and $x_l$ the reactance of the $l$-th line. Define then the $L \times L$ diagonal matrix $\mathbf{D} := \mathrm{diag}\left( \{x_l^{-1}\}_{l=1}^{N_l}\right)$; and the $L\times M$ branch-bus incidence
matrix $\mathbf{A}$, such that if its $l$-th row $\mathbf{a}_l^T$
corresponds to the branch $(m,n)$, then $[\mathbf{a}_l]_m:=+1$,
$[\mathbf{a}_l]_n:=-1$, and zero elsewhere.

Flow conservation dictates that the aggregate power injected per bus should equal the power flowing away from the bus. The DC power flow model gives rise to the \emph{nodal balance constraint}~\cite{ExpConCanBook}
\begin{equation}
\label{eq:node_bal}
\mathbf{p}_G + \mathbf{w} - \mathbf{p}_D = \mathbf{B}\bm\theta
\end{equation}
where $\bm\theta:=[\theta_1~\ldots~\theta_M]^T$ is the vector of nodal voltage phases $\{\theta_m\}_{m=1}^M$, and $\mathbf{B} := \mathbf{A}^T \mathbf{D} \mathbf{A}$ is the bus admittance matrix. Since the all-ones vector belongs to the nullspace of $\mathbf{B}$, the node balance equation \eqref{eq:node_bal} is invariant to nodal phase shifts. Hence, without loss of generality, the first bus can be the reference bus with phase set to zero, that is, $\theta_1=0$.

According again to the DC flow model, the power flows on all transmission lines can be expressed as $\mathbf{H}\bm\theta$ for $\mathbf{H} := \mathbf{D} \mathbf{A}$. Physical considerations enforce a  limit $\mathbf{f}^{\max}$ on the transmission power flows leading to the \emph{line flow constraint}
\begin{equation}
\label{eq:flow_lim}
-\mathbf{f}^{\max} \preceq \mathbf{H}\bm\theta \preceq \mathbf{f}^{\max}
\end{equation}
where $\preceq$ denotes entry-wise inequality.

Recall that the power system is dispatched several hours or even one day prior to the operation period of interest. Given a wind power generation forecast $\mathbf{z}$, the system operator wishes to schedule an injection $\mathbf{w}$ that is expected to be furnished. This requirement is captured here by allowing the vector inequality $\mathbf{z}\succeq\mathbf{w}$ to be violated with very low risk $\alpha$. Specifically, the following \emph{chance constraint} is imposed:
\begin{equation}
\label{eq:chance}
\mathrm{Prob}\left(\mathbf{z}\succeq \mathbf{w}\right) \geq 1-\alpha
\end{equation}
and typical values for the risk level $\alpha$ are 1--5$\%$.

There are two main challenges in dealing with~\eqref{eq:chance}. The first one is that the distribution of $\mathbf{z}$ is rarely known, as it is dictated by complex meteorological and harvesting technology related considerations. For the case of wind energy, simplified models for the power generated by individual wind farms are available---see e.g., \cite{LiuX10}, \cite{DhDoNAPS12}---but accounting for the spatial correlation among wind power producers still renders the distribution of $\mathbf{z}$ intractable. The second challenge is that constraint~\eqref{eq:chance} is generally nonconvex.

Let $C_m(p_{G_m})$ be the cost associated with the $m$-th thermal generator. Function $C_m(p_{G_m})$ is convex and strictly increasing, with typical forms being quadratic or piecewise linear. The scheduling problem amounts to minimizing the total production cost subject to the constraints presented earlier, that is,
\begin{subequations}
\label{eq:DCOPF-all}
\begin{align}
\min_{\mathbf{p}_G, \mathbf{w},\bm \theta}~~& \sum_{m=1}^{M} C_m(p_{G_m})
\label{eq:DCOPF-obj}\\
\text{subj.~to~~} & \mathbf{p}_G + \mathbf{w} - \mathbf{p}_D = \mathbf{B}\bm\theta
\label{eq:DCOPF-node}\\
& -\mathbf{f}^{\max} \preceq \mathbf{H}\bm\theta \preceq \mathbf{f}^{\max}
\label{eq:DCOPF-line}\\
& \mathbf{p}_{G}^{\min} \preceq \mathbf{p}_{G} \preceq \mathbf{p}_{G}^{\max}
\label{eq:DCOPF-gen}\\
& \theta_1=0
\label{eq:DCOPF-ref}\\
& \mathrm{Prob}\left(\mathbf{z}\succeq\mathbf{w}\right) \geq 1-\alpha.
\label{eq:DCOPF-chance}
\end{align}
\end{subequations}

Formulation \eqref{eq:DCOPF-all} extends to the DC optimal power flow (OPF) problem---see e.g., \cite{ChWoWa00}---to account for uncertain renewable energy injections. To this end, the scheduled renewable energy $\mathbf{w}$ is used as a basis for optimizing the power outputs of thermal generators based on~\eqref{eq:DCOPF-node}. The risk that the produced renewable energy will not be adequate to provide the scheduled one is limited as per constraint~\eqref{eq:DCOPF-chance}. If during the actual system operation the harvested renewable energy exceeds the scheduled value, then curtailment is effected.

Note that constraints~\eqref{eq:DCOPF-node}--\eqref{eq:DCOPF-gen} are linear and the objective~\eqref{eq:DCOPF-obj} is convex. Nevertheless, convexity of the overall problem~\eqref{eq:DCOPF-all} is lost due to~\eqref{eq:DCOPF-chance}. Recall also that the left-hand side of~\eqref{eq:DCOPF-chance} is difficult to be expressed as a function of the decision variable $\mathbf{w}$, while the constraint is generally nonconvex. To this end, the ensuing section develops a numerically tractable convex approximation of~\eqref{eq:DCOPF-chance}.

\section{Scenario Approximation Approach}\label{Section:SAA}

The convex approximation of~\eqref{eq:DCOPF-chance} relies on the scenario-based approximation proposed originally for robust control~\cite{Calafiore06}, and recently used for chance-constrained economic dispatch in~\cite{YuNGGG-ISGT13}. The method relies on the availability of independent samples from the distribution of $\mathbf{z}$.

Specifically, let $\{\mathbf{z}(s)\}_{s=1}^S$ denote $S$ independent samples available. The scenario approximation approach relies on substituting \eqref{eq:DCOPF-chance} with its sampled version
\begin{equation}
\label{eq:sampleconstr}
\mathbf{w} \preceq \mathbf{z}(s), \quad s=1,\ldots,S.
\end{equation}
Then, the optimization problem consisting of~\eqref{eq:DCOPF-obj}--\eqref{eq:DCOPF-gen} and~\eqref{eq:sampleconstr} is solved.

As~\eqref{eq:sampleconstr} is an approximation of~\eqref{eq:DCOPF-chance}, the question of whether the solution of the resultant optimization problem is feasible for the original problem is raised. In fact, notice that the solution of the approximate problem is a random variable, because the samples  $\{\mathbf{z}(s)\}_{s=1}^S$ are random.
Reference~\cite{Calafiore06} develops a bound on the sample size $S$ as a function of the risk level $\alpha$ which guarantees that the solution of the approximate problem is feasible for the original one with high probability.

Notice that~\eqref{eq:sampleconstr} is linear in $\mathbf{w}$, a fact that renders the overall scheduling problem convex. On the negative side, the required $S$ to achieve feasibility of the approximate solution is typically very large. This implies that the resultant optimization problem will have a very large number of constraints [cf.~\eqref{eq:sampleconstr}], which may pose significant computational burden.  It is possible to exploit the structure of the problem at hand, in order to overcome this difficulty and come up with a sample size free approximation. Specifically, it is not hard to see that~\eqref{eq:sampleconstr} is equivalent to
\begin{equation}
\label{eq:minsampleconstr}
{w}_m \leq \min_{s=1,\ldots,S}\{{z}_m(s)\}, \quad m=1,\ldots,M.
\end{equation}
A complication of this sampling mechanism is that the right-hand side of~\eqref{eq:minsampleconstr} can become very small as $S$ grows. Recall that $z_m$ is the power output of the $m$th renewable energy producer. As such, it is lower bounded by zero, and there is in fact nonzero probability that $z_m=0$. This shortcoming can drive the decision variable $w_m$ to very small values or even to zero. In a nutshell, there is a degree of conservatism inherent to the scenario approximation method.

A straightforward modification of \eqref{eq:minsampleconstr} can alleviate the aforementioned conservatism. Specifically, a small quantity $\delta_m>0$ can be added to $\min_{s=1,\ldots,S}\{{z}_m(s)\}$, in which case \eqref{eq:minsampleconstr} is surrogated by
\begin{equation}\label{eq:minlift}
{w}_m \leq \min_{s=1,\ldots,S}\{{z}_m(s)\} +\delta_m, \quad m=1,\ldots,M.
\end{equation}
The effectiveness of this adjustment will be demonstrated numerically. With $\mathbf{z}^{\text{res}} = [z_1^{\text{res}}~\ldots~z_M^{\text{res}}]$ denoting the right-hand side of~\eqref{eq:minlift}, the following problem is solved instead of~\eqref{eq:DCOPF-all}:
\vspace{.2cm}

\begin{subequations}\label{eq:apprDCOPF-all}
\hspace{-.5cm}
\fbox{
 \addtolength{\linewidth}{-2\fboxsep}%
 \addtolength{\linewidth}{-2\fboxrule}%
 \begin{minipage}{0.98\linewidth}
\begin{align}
\min_{\mathbf{p}_G, \mathbf{w},\bm \theta}~~& \sum_{m=1}^{M} C_m(p_{G_m})
\label{eq:apprDCOPF-obj}\\
\text{subj.~to~~} & \mathbf{p}_G + \mathbf{w} - \mathbf{p}_D = \mathbf{B}\bm\theta
\label{eq:apprDCOPF-node}\\
& -\mathbf{f}^{\max} \preceq \mathbf{H}\bm\theta \preceq \mathbf{f}^{\max}
\label{eq:apprDCOPF-line}\\
& \mathbf{p}_{G}^{\min} \preceq \mathbf{p}_{G} \preceq \mathbf{p}_{G}^{\max}
\label{eq:apprDCOPF-gen}\\
& \theta_1=0
\label{eq:apprDCOPF-ref}\\
& \mathbf{w} \preceq \mathbf{z}^{\text{res}}.
\label{eq:apprDCOPF-scen}
\end{align}
\end{minipage}
}
\end{subequations}

\vspace{.2cm}

\section{Numerical Tests}\label{sec:Numericalresults}
\begin{figure*}[t]
\centering
\includegraphics[width=0.7\linewidth]{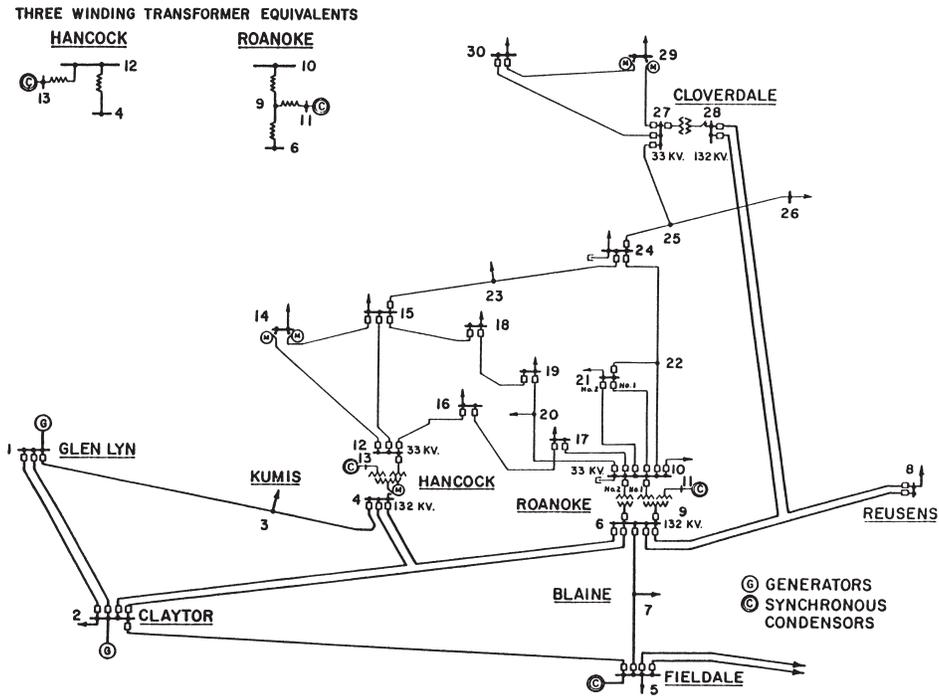}
\caption{Seven wind farms have been added to the IEEE 30-bus grid system~\cite{PSTCA}.}
\label{fig:ieee30}
\end{figure*}

The performance of the novel scheduling approach is corroborated via numerical tests using the IEEE 30-bus power system~\cite{PSTCA}. The latter includes 41 transmission lines and 6 conventional generators residing at buses $\{1,2,13,22,23,27\}$ (cf.~Fig.\ref{fig:ieee30}). Load demands, guadratic generation costs, generator capacities, and transmission line ratings, are all specified in \cite{MATPOWER}. Seven wind farms have been added on buses $\{1,2,5,9,15,24,30\}$. The convex problem in \eqref{eq:apprDCOPF-all} is solved using the \texttt{CVX} package and the \texttt{SDPT3} solver \cite{cvx}, \cite{sdpt3}.

To simulate wind farm operation, real data originally provided for a wind energy forecasting competition organized by Kaggle platform were utilized~\cite{kaggle}. Among other data, the specific dataset contains the actual hourly power output of seven wind farms over three years. To eliminate possible non-stationarities, only the interval from May $1^{\text{st}}$ to June $26^{\text{th}}$ of 2012 was considered, yielding a total of 589 hours due to missing entries.

Wind power outputs have been normalized per farm due to privacy concerns. To preserve the total installed generation capacity fixed after adding the wind farms, the conventional capacity is scaled down by 80\%. Then, all wind farm outputs are scaled to contribute equally to the rest of the installed capacity, hence yielding a 20\% wind energy penetration.

Recall that the developed scenario approximation-based scheduling requires drawing independent samples from the wind energy forecast $\mathbf{z}$. As a proof of concept, it is assumed here that $\mathbf{z}$ is Gaussian distributed. Its expected value is considered to be the actual wind power generated. A ``low-wind'' and a ``high-wind'' scenario were considered. The low-wind scenario yields $\boldsymbol{\mu}_l= [1.15~1.37~0.47~1.05~1.45~1.64~0.00]^T$ and corresponds to May 19$^{\text{th}}$ at 8~a.m. The high-wind scenario has $\boldsymbol{\mu}_h = [6.00~0.31~7.66~8.01~8.42~ 8.44~8.46]^T$ and is observed on May 22$^{\text{nd}}$ at 8~a.m. To model correlation across farms, it is further postulated that the covariance of $\mathbf{z}$ is that of the wind farm power outputs. The latter is empirically estimated as the sample covariance and it is denoted by $\hat{\mathbf{\Sigma}}$. Samples of $\mathbf{z}$ can then be drawn from $\mathcal{N}(\boldsymbol{\mu}_l,\hat{\mathbf{\Sigma}})$ and $\mathcal{N}(\boldsymbol{\mu}_h,\hat{\mathbf{\Sigma}})$, respectively for the two scenarios.

Before solving \eqref{eq:apprDCOPF-all}, the boosting parameters $\{\delta_m\}_{m=1}^M$ introduced in \eqref{eq:minlift} must be selected. An intuitive and easily-implementable heuristic for doing so is described next. Instead of constraining $w_m$ to be no larger than \emph{all} samples $z_m(s)$ as dictated by \eqref{eq:minsampleconstr}, it is natural to require $w_m$ to be no larger than only the $(1-\alpha)\%$ largest samples. Algorithmically, if $\left\{\{z_m^{[s]}\}: z_m^{[1]}\geq z_m^{[2]} \geq \cdots \geq z_m^{[S]}\right\}$ denote the order statistics of the original samples $\{z_m(s)\}_{s=1}^{S}$ for $m=1,\ldots,M$, the right-hand side of \eqref{eq:minlift} can be selected as ${z}_m^{\text{res}} = z_m^{[\lceil(1-\alpha)\times S\rceil]}$. Negative-valued entries of $\mathbf{z}^{\text{res}}$ are truncated to zero.

Dispatching the IEEE 30-bus power system for a risk level of $\alpha=0.05$ yields the optimal costs listed in Table~\ref{tab:bad-good}. The Lagrange multipliers corresponding to \eqref{eq:apprDCOPF-node}, also known as \emph{locational marginal prices} (LMPs), are also listed in the same table. LMPs are important components of electricity markets since they represent the cost of selling or buying electricity at a particular bus; see e.g., \cite{KirschenStrbac}, \cite{ExpConCanBook}. Due to lack of transmission line congestion, all LMPs turn out to be equal to the value provided in Table~\ref{tab:bad-good}. The high-wind scenario attains lower cost and LMPs than the low-wind scenario, since less conventional power is needed when more free wind power is available. It is worth mentioning that due to the risk-aware constraint, the low-wind scenario essentially boils down to scheduling with no wind power at all.

\begin{table}[t]
\centering
\caption{Optimal costs and LMPs for high-/low-wind scenarios $(\alpha=0.05)$.}\label{tab:bad-good}
    \begin{tabular}{ c || c | c }
    \backslashbox{Scenario}{}
                        &Cost        &LMP \\  \hline

    High-wind         &481.42    &364.97    \\  \hline
    Low-wind        &565.21    &378.91   \\
\hline
    \end{tabular}
\end{table}

\begin{figure}[t]
\centering
\includegraphics[width=0.4\textwidth]{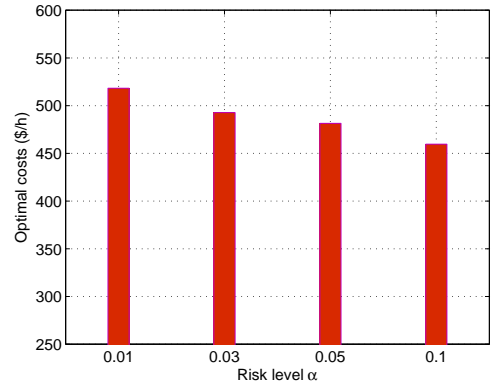}
\caption{Optimal costs for varying $\alpha$ (high-wind scenario).}
\label{fig:cost_Wgood}
\end{figure}
\begin{figure}[t]
\centering
\includegraphics[width=0.4\textwidth]{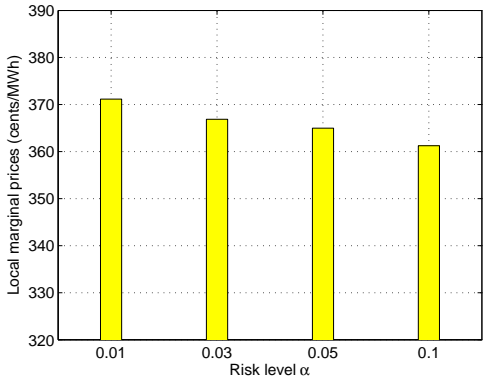}
\caption{Locational marginal prices for varying $\alpha$ (high-wind scenario).}
\label{fig:LMP_Wgood}
\end{figure}

\begin{table}[t]
\centering
\caption{Prescribed risk level and actual risk (high-wind scenario).}\label{tab:true-risk}
    \begin{tabular}{ c || c | c | c | c }
    \hline
    $\alpha$          &0.01    &0.03    &0.05    &0.1  \\  \hline
    Actual risk          &0.0072    &0.0075    &0.0076    &0.0087    \\
\hline
    \end{tabular}
\end{table}

Figs.~\ref{fig:cost_Wgood} and \ref{fig:LMP_Wgood} illustrate the effect of the prescribed risk level $\alpha$ on the optimal costs and the LMPs, respectively. The optimal net cost decreases with increasing $\alpha$, since higher risk allows more wind power to be committed.

To justify the heuristic boosting procedure, the risk incurred by the $\mathbf{w}$ minimizing \eqref{eq:apprDCOPF-all} is empirically evaluated by drawing $10^5$ independent wind forecast samples $\mathbf{z}$, and checking whether \eqref{eq:chance} holds. Table~\ref{tab:true-risk} shows the validation results. The actual risk is always smaller than the predefined one, hence numerically validating the boosting step.

\begin{table}[t]
\centering
\caption{Optimal costs for varying $\alpha$ and $\beta$ (high-wind scenario).}
\label{tab:albe}
    \begin{tabular}{|c||*{4}{c|}}\hline
    \backslashbox{$\alpha$}{$\beta$}
                         &1.05 &1.1  &1.2 &1.3 \\  \hline
    0.01                    &553.2936  &589.0077  &662.2263  &738.4597 \\
    0.03                     &527.7588  &565.2449  &636.9299  &712.1263 \\
    0.05                     &515.6234  &549.3798  &623.0269  &697.1422 \\
    0.1                     &496.8349  &530.5984  &603.5701  &677.1949  \\
\hline
    \end{tabular}
\end{table}

\begin{figure}[t]
\centering
\includegraphics[width=0.45\textwidth]{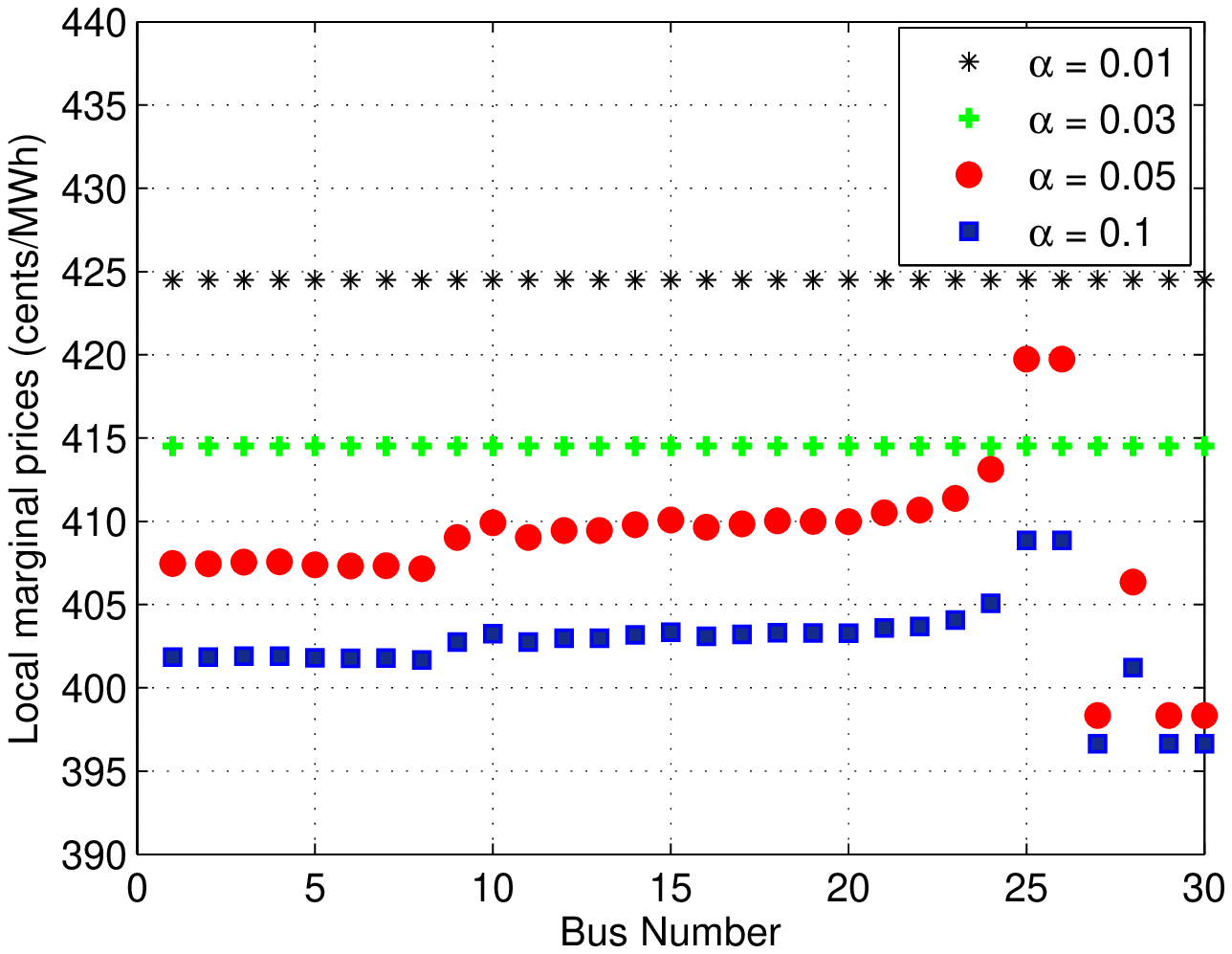}
\vspace*{-1em}
\caption{Locational marginal prices for $\beta=1.330$.}
\label{fig:LMP-congest}
\end{figure}

\begin{figure}[t]
\centering
\includegraphics[width=0.45\textwidth]{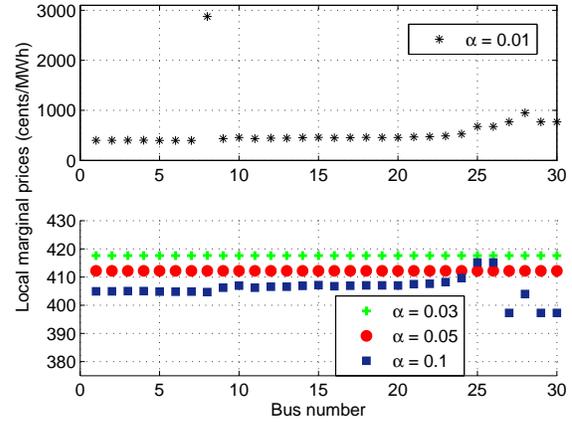}
\vspace*{-1em}
\caption{Locational marginal prices for $\beta=1.342$.}
\label{fig:LMP-congest2}
\end{figure}

The effect of the risk level $\alpha$ on LMPs under transmission network congestion is investigated next. To simulate congestion, load demand at all buses is scaled up by $\beta$. The optimal costs listed in Table~\ref{tab:albe} decrease with decreasing $\beta$ and/or increasing $\alpha$, as expected. The corresponding 30 LMPs (one per bus) obtained when $\beta=1.330$ and $\beta=1.342$ and for varying values of $\alpha$ are plotted in Figs.~\ref{fig:LMP-congest} and~\ref{fig:LMP-congest2}, respectively. The two figures indicate that high risk levels result in lower prices in general. However, by varying $\alpha$ and $\beta$, different congestion patterns may occur due to the grid topology.

\section{Conclusions}\label{sec:Conclusion}
Network-constrained economic dispatch with multiple wind farms was considered in this paper. A risk-constrained optimization problem was formulated based on the loss-of-load probability over all wind farm injection points. To address the imperfect knowledge of the wind power forecasts, a scenario approximation technique via Monte Carlo sampling was proposed. The attractive features and practical impact of this work are two-fold: i) the scenario approach enables economic and risk-limited scheduling of
smart grids with increasingly higher renewable energy penetration, without relying on specific probabilistic assumptions about the renewable generation; and, ii) the special problem structure renders the approach applicable to large-scale problems. Multi-period formulations constitute an interesting future direction.

\bibliographystyle{IEEEtran}
\bibliography{biblio}

\begin{thebibliography}{10}
\providecommand{\url}[1]{#1}
\csname url@samestyle\endcsname
\providecommand{\newblock}{\relax}
\providecommand{\bibinfo}[2]{#2}
\providecommand{\BIBentrySTDinterwordspacing}{\spaceskip=0pt\relax}
\providecommand{\BIBentryALTinterwordstretchfactor}{4}
\providecommand{\BIBentryALTinterwordspacing}{\spaceskip=\fontdimen2\font plus
\BIBentryALTinterwordstretchfactor\fontdimen3\font minus
  \fontdimen4\font\relax}
\providecommand{\BIBforeignlanguage}[2]{{%
\expandafter\ifx\csname l@#1\endcsname\relax
\typeout{** WARNING: IEEEtran.bst: No hyphenation pattern has been}%
\typeout{** loaded for the language `#1'. Using the pattern for}%
\typeout{** the default language instead.}%
\else
\language=\csname l@#1\endcsname
\fi
#2}}
\providecommand{\BIBdecl}{\relax}
\BIBdecl

\bibitem{GWEC}
{GWEC}, ``Global wind statistics 2012,'' Feb. 2013, [Online]. Available:
  \url{http://www.gwec.net/wp-content/uploads/2013/02/GWEC-PRstats-2012_english.pdf}.

\bibitem{DOE08}
``20\% wind energy by 2030: Increasing wind energy's contribution to {U.S.}
  electricity supply,'' July 2008, [Online]. Available:
  \url{http://www1.eere.energy.gov/wind/pdfs/41869.pdf}.

\bibitem{LiuX10}
X.~Liu and W.~Xu, ``Economic load dispatch constrained by wind power
  availability: A here-and-now approach,'' \emph{IEEE Trans. on Sustainable
  Energy}, vol.~1, no.~1, pp. 2--9, Apr. 2010.

\bibitem{HetzerYB08}
J.~Hetzer, C.~Yu, and K.~Bhattarai, ``An economic dispatch model incorporating
  wind power,'' \emph{IEEE Trans. on Energy Conver.}, vol.~23, no.~2, pp.
  603--611, Jun. 2008.

\bibitem{XieGZG11}
L.~Xie, Y.~Gu, X.~Zhu, and M.~G. Genton, ``Power system economic dispatch with
  spatio-temporal wind forecasts,'' in \emph{Proc. of IEEE EnergyTech},
  Cleveland, OH, May 2011.

\bibitem{YuNGGG-ISGT13}
Y.~Zhang, N.~Gatsis, and G.~B. Giannakis, ``Risk-constrained energy management
  with multiple wind farms,'' in \emph{Proc. Innovative Smart Grid
  Technologies}, Washington, D.C., Feb. 2013.

\bibitem{ChWoWa00}
R.~D. Christie, B.~F. Wollenberg, and I.~Wangensteen, ``Transmission management
  in the deregulated environment,'' \emph{Proc. of the IEEE}, vol.~88, no.~2,
  pp. 170--195, Feb. 2000.

\bibitem{BiChHa12}
D.~Bienstock, M.~Chertkov, and S.~Harnett, ``Chance constrained optimal power
  flow: Risk-aware network control under uncertainty,'' Sept. 2012, [Online].
  Avaialble: \url{http://arxiv.org/pdf/1209.5779.pdf}.

\bibitem{SjGaTo12}
E.~Sj\"{o}din, D.~F. Gayme, and U.~Topcu, ``Risk-mitigated optimal power flow
  for wind powered grids,'' in \emph{Proc. American Control Conf.},
  Montr\'{e}al, Canada, June 2012, pp. 4431--4437.

\bibitem{kaggle}
\BIBentryALTinterwordspacing
Global energy forecasting competition 2012 - wind forecasting. [Online].
  Available: \url{http://www.kaggle.com/c/GEF2012-wind-forecasting}
\BIBentrySTDinterwordspacing

\bibitem{PSTCA}
\BIBentryALTinterwordspacing
Power systems test case archive. Univ. of Washington. [Online]. Available:
  \url{http://www.ee.washington.edu/research/pstca/}
\BIBentrySTDinterwordspacing

\bibitem{MATPOWER}
R.~D. Zimmerman, C.~E. Murillo-Sanchez, and R.~J. Thomas, ``{MATPOWER}:
  steady-state operations, planning and analysis tools for power systems
  research and education,'' \emph{IEEE Trans. on Power Syst.}, vol.~26, no.~1,
  pp. 12--19, Feb. 2011.

\bibitem{ExpConCanBook}
A.~G\'{o}mez-Exp\'{o}sito, A.~J. Conejo, and C.~Canizares, Eds., \emph{Electric
  Energy Systems, Analysis and Operation}.\hskip 1em plus 0.5em minus
  0.4em\relax Boca Raton, FL: CRC Press, 2009.

\bibitem{DhDoNAPS12}
S.~V. Dhople and A.~D. Dominguez-Garcia, ``A framework to determine the
  probability density function for the output power of wind farms,'' in
  \emph{Proc. North American Power Symposium}, Urbana, IL, Sept. 2012.

\bibitem{Calafiore06}
G.~Calafiore and M.~Campi, ``The scenario approach to robust control design,''
  \emph{IEEE Trans. Automat. Contr.}, vol.~51, pp. 742--753, 2006.

\bibitem{cvx}
{CVX Research Inc.}, ``{CVX}: Matlab software for disciplined convex
  programming, version 2.0 (beta),'' \url{http://cvxr.com/cvx}, Sep. 2012.

\bibitem{sdpt3}
SDPT3, \url{http://www.math.nus.edu.sg/~mattohkc/sdpt3.html}.

\bibitem{KirschenStrbac}
D.~Kirschen and G.~Strbac, \emph{Power System Economics}.\hskip 1em plus 0.5em
  minus 0.4em\relax West Sussex, England: Wiley, 2010.

\end{thebibliography}
\end{document}